%October 13, 2000
\magnification 1000
\input amstex
\documentstyle{amsppt}
\vsize 8.25in
\voffset 1cm
\hoffset 1cm
\topmatter
\rightheadtext{On the convergence of formal CR maps between hypersurfaces}
\leftheadtext{Jo\"el Merker}
\title Convergence of formal invertible CR mappings
between minimal holomorphically nondegenerate
real analytic hypersurfaces\endtitle
\author Jo\"el Merker
\endauthor
\address Laboratoire d'Analyse, Topologie et Probabilit\'es,
Centre de Math\'ematiques et d'Informatique, UMR 6632,
39 rue Joliot Curie,
F-13453 Marseille Cedex 13, France.
Fax: 00 33 (0)4 91 11 35 52\endaddress
\email merker\@cmi.univ-mrs.fr\endemail
\thanks
\endthanks

\loadeufm

\def\N{{\Bbb N}}

\def\C{{\Bbb C}}

\define \dl{[\![}
\define \dr{]\!]}
\def\dim{\hbox{dim}}
\def\sumb{\sum_{\beta\in \N_*^{n-1}}}
\def\sumg{\sum_{\gamma\in \N_*^{n-1}}}
\def\n*{\N_*^{n-1}}
\def\1{{\text{\bf 1}}}
\def\v{\vert}
\def\n{\vert\vert}

%\abstract

%\centerline{\bf Table of contents~:}

%\smallskip

%{\bf \S1.~Introduction and statement of the results \dotfill 1.}

%{\bf \S2.~Preliminaries and notations \dotfill 4.}

%{\bf \S3.~Minimality and holomorphic nondegeneracy\dotfill 7.}

%{\bf \S4.~Formal versus analytic \dotfill 8.}

%{\bf \S5.~Classical reflection identities \dotfill 9.}

%{\bf \S6.~Convergence of the reflection function on $\Cal S_0^1$ \dotfill 12.}

%{\bf \S7.~Convergence of the jets of the reflection function on $\Cal S_0^1$ \dotfill 14.}

%{\bf \S8.~Convergence of the formal mapping \dotfill 17.}

%\endabstract

\endtopmatter

\document

\head \S 1. Introduction and statement of the results \endhead

\subhead 1.1.~Main theorem\endsubhead
We establish here the following assertion (here and in what follows,
the symbol $\Cal F$ relates to objects and maps of the 
{\it formal category}\,:

\proclaim{Theorem~1.2}
Any invertible $(${\it i.e.}~with nonzero 
Jacobian determinant at $p)$ formal CR mapping $h\: (M,p)\to_{\Cal F}
(M',p')$ between two germs of minimal real analytic
hypersurfaces in $\C^n$, $n\geq 2$, is convergent \text{\rm if and
only if} $(M',p')$ is holomorphically nondegenerate.
\endproclaim

\noindent
(The reader is referred to the monograph [3] and to the articles
[2,4,10,12] for background material). This theorem provides a necessary and
sufficient condition for the convergence of an invertible formal CR
map of hypersurfaces. The necessity appears in a natural way ({\it
see} Proposition~1.5 below). Geometrically, holomorphic nondegeneracy has
clear geometric signification\,: it means that there exist {\it no}
holomorphic tangent vector field to $(M',p')$ and it is equivalent to the
{\it nonexistence} of a local complex analytic foliation of $(\C^n, p')$
tangent to $(M', p')$. As matters stand, such a kind of
characterization for the regularity of CR maps happens to be known
already, but only in case where at least {\it one} of the two
hypersurfaces is algebraic, see {\it e.g.}  [5,6,13] (in fact, 
in the algebraic case, one can study the problem by means of 
classical ``polynomial identities'' in
the spirit of Baouendi-Jacobowitz-Treves).  But it was known that the
truly real analytic case involves deeper investigations.

\subhead 1.3.~Brief history\endsubhead
Formal invertible CR mappings $h\: (M,p)\to_{\Cal F} (M', p')$ between
two germs of real analytic hypersurfaces in $\C^n$ have been proved to
be convergent in many circumstances. Firstly, in 1974 by Chern-Moser, 
assuming that $(M',p')$ is Levi-nondegenerate.  Secondly, in
1997 by Baouendi-Ebenfelt-Rothschild in the important article [2],
assuming that $h$ is invertible ({\it i.e.} with nonzero Jacobian at
$p$) and that $(M', p')$ is finitely nondegenerate at $p'$. And more
recently in 1999, by Baouendi-Ebenfelt-Rothschild [4], assuming
for instance (but this work also contains other results) that
$(M',p')$ is essentially finite, that $(M,p)$ is minimal and that $h$
is not totally degenerate, a result which is valid in arbitrary
codimension. (Again, the reader may consult [3] for essential
background on the subject, for definitions, concepts and tools and
also [10] for related topics.) In summary, the above-mentioned results
have all exhibited various sufficient conditions.

\remark{Remark}
After a first version of this preprint was finished and distributed,
the author received a preprint (now published) [12] where a statement
similar to Theorem~1.2 was proved.  Further, after dropping the
assumption of holomorphic nondegeneracy on $(M',p')$, the convergence
of the reflection function is established~in~it.
\endremark

\subhead 1.4.~Necessity\endsubhead
On the other hand, it is known (essentially since 1995, {\it cf.}~[5])
that holomorphic nondegeneracy of the hypersurface $(M',p')$
constitutes a {\it natural necessary condition} for $h$ to be
convergent, according to an important observation due to
Baouendi-Rothschild [2,3,5] (this observation followed naturally from
the characterization by Stanton of the finite-dimensionality of the
space of infinitesimal CR automorphisms of $(M,p)$ [16]~; Stanton's
discovery is fundamental in the subject).  We may restate this
observation as follows ({\it see} its proof in the end of \S4).

\proclaim{Proposition~1.5}
If $(M', p')$ is holomorphically degenerate, then there exists a
nonconvergent formal invertible CR self map of $(M',p')$,
which is simply of the form $\C^n\ni t'\mapsto_{\Cal F} \exp (\varpi'(t')
L')(t')\in \C^n$, where $L'$ is a nonzero holomorphic tangent vector
to $(M',p')$ and where $\varpi'(t')\in \C\dl t'\dr$, $\varpi'(0)=0$, is
nonconvergent.
\endproclaim

\remark{Remark}
A geometric way to interpret this nonconvergent map would be to say
that it flows in nonconvergent complex time along the complex analytic
foliation induced by $L'$, which is tangent to $M'$. Similar
obstructions for the algebraic mapping problem stem from the existence of
complex analytic (or algebraic) foliations tangent to $(M', p')$, see
{\it e.g.} [5,6]. Again, this shows that the geometric notion of 
holomorphic nondegeneracy discovered by Stanton is crucial in the field. 
\endremark

\subhead 1.6.~Jets of Segre varieties\endsubhead
The holomorphically nondegenerate hypersurfaces are considerably more
general and more complicated to handle than Levi-nondege\-ne\-rate
ones [14,15,17], finitely nondegenerate ones [2], essentially finite
ones [3,4] or even Segre nondegenerate ones [10]. The explanation
becomes simple after a reinterpretation of these conditions in the
spirit of the important geometric definition of jets of Segre
varieties due to Diederich-Webster [7]. In fact, these five distinct
nondegeneracy conditions manifest themselves directly as nondegeneracy
conditions of the morphism of $k$-th jets of Segre varieties attached
to $M'$, which is an invariant {\it holomorphic} map defined on its
extrinsic complexification $\Cal M'=(M')^c$ (we follow the notations
of \S2). In local holomorphic normal coordinates $t'=(w',z')\in
\C^{n-1}\times \C$, vanishing at $p'$ with $\tau':=(\zeta',\xi')\in
\C^{n-1}\times \C$ denoting the complexied coordinates $(w',z')^c$,
such that the holomorphic equation of the extrinsic complexification
$\Cal M'$ is written $\xi'=z'-i\Theta'(\zeta',t')=z'-i\sum_{\gamma\in
\N_*^{n-1}} {\zeta'}^{\gamma} \, \Theta_\gamma'(t')$ 
({\it cf.}~\thetag{2.4}), the conjugate
complexified Segre variety is defined by $\underline{\Cal
S}_{t'}':=\{\tau'\: \xi'=z'-i\Theta'(\zeta',t')\}$ (here, $t'$ is
fixed\,; {\it see} [9] for a complete exposition of the geometry of
complexified Segre varieties) and the jet of order $k$ of the
complex $(n-1)$-dimensional manifold
$\underline{S}_{t'}'$ at the point $\tau'\in\underline{S}_{t'}'$ 
defines a holomorphic map
$$
\varphi_k'\: {\Cal M}' \ni (t', \tau') \mapsto j_{\tau'}^k
\underline{S}_{t'}' \in
\C^{n+N_{n-1,k}}, \ \ \ \ \ N_{n-1,k}= {(n-1+k)! \over (n-1)! \, k!},
\tag 1.7
$$ 
given explicitely in terms of such a defining equation by a collection
of power series\,:
$$
\varphi_k'(t',\tau'):=
j_{\tau'}^k \underline{S}_{t'}'=(\tau',\{
\partial_{\zeta'}^\beta\left[\xi'-
z'+i\Theta'(\zeta',t')\right]\}_{\beta\in\N^{n-1},\v\beta\v\leq k}).
\tag 1.8
$$ 
For $k$ large enough, the various possible properties of this
holomorphic map govern some different ``{\it nondegeneracy conditions}''
on $M'$ which are appropriate for some generalizations of the
Lewy-Pinchuk reflection principle. Let ${p'}^c:=(p',\bar p')\in
\Cal M'$. We give here an account of five
conditions, which can be understood as definitions\,:
$$
\aligned
& 
(M',p') \ \hbox{is {it Levi-nondegenerate}} \ \hbox{at} \ p' 
\ \ \ \ \ \ \ \ \ \ \ \ \ \ \ \ \ \ \ \ \ \ \ \ \\
&
\Longleftrightarrow \ \varphi_1' \ \hbox{is an immersion at} \ 
{p'}^c.\\
\endaligned
\tag I
$$\vskip -0.2cm
$$
\aligned
& 
(M',p') \ \hbox{is {\it finitely nondegenerate}} \ \hbox{at} \ p' \\
& 
\Longleftrightarrow \ \exists \ k_0\in \N_*, \
\varphi_k' \ \hbox{is an immersion at} \ 
{p'}^c, \ \forall\ k\geq k_0.
\endaligned
\tag II
$$\vskip -0.2cm
$$
\aligned
& 
(M',p') \ \hbox{is {\it essentially finite}} \ \hbox{at} \ p' \\
& 
\Longleftrightarrow \ \exists \ k_0\in \N_*, \
\varphi_k' \ \hbox{is a finite holomorphic map at} \ 
{p'}^c, \ \forall\ k\geq k_0.
\endaligned
\tag III
$$\vskip -0.2cm
$$
\aligned
& 
(M',p') \ \hbox{is {\it S-nondegenerate}} \ \hbox{at} \ p' \\
& 
\Longleftrightarrow \ \exists \ k_0\in \N_*, \
\varphi_k'|_{{\Cal S}_{\bar p'}} \ \hbox{is of generic rank} \ 
\hbox{dim}_{\C} {\Cal S}_{\bar p'}=n-1, \ \forall\ k\geq k_0.
\endaligned
\tag IV
$$\vskip -0.2cm
$$\aligned
& 
(M',p') \ \hbox{is {\it holomorphically nondegenerate}} \ \hbox{at} \ p' \\
& 
\Longleftrightarrow \ \exists \ k_0\in \N_*, \
\varphi_k' \ \hbox{is of generic rank} \ 
\hbox{dim}_{\C} {\Cal M}'=2n-1, \ \forall\ k\geq k_0.
\endaligned
\tag V
$$

\remark{Remarks}~{\bf 1.}~It follows from the biholomorphic
invariance of Segre varieties that two Segre morphisms of $k$-jets
associated to two different local coordinates for $(M',p')$ are 
intertwined by a local biholomorphic map of $\C^{n+N_{n-1,k}}$.
Consequently, the properties of $\varphi_k'$ are invariant.

{\bf 2.}~The condition (I) is classical.  The condition (II) is
studied by Baouen\-di-Ebenfelt-Rothschild [2,3] and appeared already
in Pinchuk's thesis, in Diederich-Webster [7] and in some of Han's
works.  The condition (III) appears in Diederich-Webster [7] and was
studied by Baouendi-Jacobowitz-Treves and by Diederich-Fornaess.  The
condition (IV) seems to be new and appears in [10]. The condition (V)
was discovered by Stanton in her concrete study of infinitesimal CR
automorphisms of real analytic hypersurfaces ({\it see} [16] and the
references therein) and is equivalent to the {\it nonexistence} of a
holomorphic vector field with {\it holomorphic} coefficients tangent
to $(M',p')$.  We claim that it is easy to show that $(I) \
\Rightarrow \ (II) \ \Rightarrow \ (III)
\ \Rightarrow \ (IV) \ \Rightarrow \ (V)$ (only the implication 
$(IV) \ \Rightarrow \ (V)$ is not straightforward, {\it see}
Lemma~5.15 below for a proof). Finally, this stratification is the
same, word by word, in higher codimension.

\endremark

\subhead 1.9.~A general commentary \endsubhead  
To confirm evidence of the strong differences
between these five levels of nondegeneracy, let us point out a very
clear fact\,: the immersive or finite local holomorphic maps
$\varphi\: (X,p)\to (Y, q)$ between complex manifolds with $\dim_{\C}
X\leq\dim_\C Y$ are very rare in the set of maps of generic rank equal
to $\dim_\C X$, or even in the set of maps having maximal generic rank
$m$ over a submanifold $(Z,p)\subset (X,p)$ of positive dimension
$m\geq 1$. Thus condition (V) is by far the most general. Furthermore,
an important difference between (V) and the other conditions is that
{\it $(V)$ is the only condition which is nonlocal}, in the sense that
it happens to be satisfied at every point if it is satisfied at a
single point only, provided, of course, that $(M',p')$ is
connected. On the contrary, it is obvious that the other four
conditions are really local\,: even though they happen to be satisfied
at one point, there exist in general many other points where they fail
to be satisfied. In this concern, let us recall that any $(M',p')$
satisfying (V) must satisfy (II) locally -- hence also (III) and (IV)
-- over a Zariski dense open subset of points of $(M',p')$ (this
important fact is proved in [3]).  Therefore, the points satisfying
(III) but not (II), or (IV) but not (III), or (V) but not (IV), can
appear to be more and more exceptional and rare from the point of view
of a point moving at random in $(M',p')$, but however, from the point
of view of local analytic geometry, which is the adequate viewpoint in
this matter, they are more and more generic and general, in truth.

\remark{Remark}~An 
important feature of the theory of CR manifolds is to propagate the
properties of CR functions and CR maps along Segre chains, when
$(M,p)$ is minimal, like iteration of jets [3], support of CR
functions, {\it etc.} Based on this heuristic idea, and believing that
the generic rank of the Segre morphism over a Segre variety is a
propagating property, I have claimed in February 1999
that any real analytic $(M',p')$ which is minimal at $p'$ happens
to be holomorphically nondegenerate if and only if it is Segre
nondegenerate at $p'$. This is not true for a general $(M',p')$ as is
shown for instance by an example from [4]\,: we take in $\C^3$
equipped with affine coordinates $(z_1',z_2',z_3')$ ({\it see} Lemmas~3.3
and~5.15 for a checking)\,:
$$
M' : \ \ 
y_3'=\vert z_1' \vert^2 \vert 1+ z_1' \bar{z}_2'\vert^2 
(1+\hbox{Re} (z_1'\bar{z}_2'))^{-1} -x_3'\,
\hbox{Im} (z_1'\bar{z}_2') 
(1 +\hbox{Re} (z_1'\bar{z}_2'))^{-1}.
\tag 1.10
$$
\endremark

\subhead 1.11.~Summary of the proof\endsubhead 
To the mapping $h$, we will associate the so-called invariant {\it
reflection function} ${\Cal R}_h'(t, \bar\nu')$ as a $\C$-valued map
of $(t,\bar\nu')\in (\C^n, p) \times (\overline\C^n, \bar{p})$ which
is a series {\it a priori} only formal in $t$ and holomorphic in
$\bar\nu'$ (the interest of studying the reflection function {\it
without any nondegeneracy condition on $(M',p')$} has been pointed out
for the first time by the author and Meylan in [11]). We prove in a
first step that $\Cal R_h'$ and all its jets with respect to $t$
converge on the first Segre chain. Then using Artin's approximation
theorem [1] (the interest of this theorem of Artin for the subject has
been pointed out by Derridj in 1985) and holomorphic nondegeneracy of
$(M', p')$, we establish that the formal CR map $h$ converges on the
second Segre chain. Finally, the minimality of $(M,p)$ together with a
theorem of Gabrielov reproved elementarily by Eakin and Harris [8] will both
imply that $h$ is convergent in a neighborhood of $p$.

\subhead 1.12.~Closing remark \endsubhead
The first version of this preprint was achieved in December 1999 and
then distributed to some specialists in the field in January 2000.
The author decided not to submit the present work, essentially because
he was seeking the same result in arbitrary codimension (afterwards,
in the beginning of March 2000 he received the prepublication of the
article [12], written meanwhile). Per chance, it appeared soon to the
author that the main ideas of proof in \S8 below were helpful to get a
much stronger result. Thus, in a preprint achieved in late May 2000
and entitled {\it \'Etude de la convergence de l'application de
sym\'etrie formelle} (in French\,; {\tt arXiv.org/abs/math/0005290}), the
author actually proves that the formal reflection mapping associated
with a formal CR equivalence between two arbitrary minimal CR-generic
manifolds in $\C^n$ of arbitrary codimension is convergent (it remains
thus now to study the nonminimal case). In fact, it appears that a
very slight modification of the techniques of the present paper yields
this stronger result naturally\,: it suffices to make a direct
application of the Artin approximation Theorem~4.2 below to
eqs~\thetag{8.5} with formal solutions $h(w,i\bar \Theta(w,\zeta,0))$
and then to proceed by induction, which we {\it do not do} here but do
in the paper of May 2000. Thus, at the present time, the author
believes that his present partial result of December 1999 has enough
interest in itself to be made public. Further, the preprint of May
borrows directly some background material to the present paper. To end
up these comments, the author would like to mention that the
inductional study of the jets of the formal CR map (or of the
reflection mapping) that he strongly uses in his works was discovered
by Baouendi-Ebenfelt-Rothschild in their study of the algebraicity of
holomorphic mappings [3] and then transferred by them to formal CR
maps in the fundamental article [2].  The known nontrivial open
problem in the subject was to treat the holomorphically nondegenerate
hypersurfaces to get an optimal sufficient condition of convergence.

\head \S2. Preliminaries and notations \endhead

\subhead 2.1.~Defining equations\endsubhead
We shall identify any germ of a real analytic hypersurface, with a
small representative of it. Thus, we shall assume constantly that we
are given two small local real analytic manifold-pieces $(M,p)$ and
$(M',p')$ of hypersurfaces in $\C^n$ with centered points $p\in M$ and
$p'\in M'$.  We first choose local holomorphic coordinates $t=(w,z)\in
\C^{n-1}\times \C$, $z=x+iy$ and $t'=(w',z') \in \C^{n-1} \times \C$,
$z'=x'+iy'$, vanishing at $p$ and at $p'$ such that the tangent spaces
to $M$ and to $M'$ at $0$ are given by $\{y=0\}$ and by $\{y'=0\}$ in
these coordinates. By this choice, we carry out ({\it cf.} [3]) the
equations of $M$ and of $M'$ under the form
$$
M\: \ \ \
z=\bar z + i\bar\Theta(w, \bar w,\bar z) \ \ \ 
\hbox{and} \ \ \ M'\: \ \ \ 
z'=\bar z' + i\bar\Theta'(w', \bar w',\bar z'), \ \ \
\tag 2.2
$$
where the power series $\bar\Theta$ and $\bar\Theta'$ converge
normally in $(2r
\Delta)^{2n-1}$ for some small $r>0$. We denote by
$\v t \v : = \sup_{1\leq i \leq n} \vert t_i \vert $ the polydisc
norm, so that $(2r\Delta)^{2n-1}= \{\v w \v,
\v \bar w\v, \v \bar z\v < 2r \}$. Here, if we denote by 
$\tau:=(\bar{t})^c:=(\zeta,\xi)$ the extrinsinc complexification of
the variable $\bar{t}$, the equations of the complexified
hypersurfaces $\Cal M:=M^c$ and $\Cal M':=(M')^c$ are simply obtained
by complexifying the eqs.~\thetag{2.2}\,:
$$
{\Cal M} \: \ \ \ 
z=\xi+i\bar\Theta(w,\zeta,\xi) \ \ \ \hbox{and} \ \ \ 
{\Cal M}' \: \ \ \ 
z'=\xi'+i\bar\Theta'(w',\zeta',\xi').
\tag 2.3
$$
As in [3], we shall assume for convenience that the coordinates
$(w,z)$ and $(w',z')$ are {\it normal}, {\it i.e.} that they are
already straightened in order that $\Theta(\zeta,0, z)\equiv 0$,
$\Theta(0,w,z)\equiv 0$ and $\Theta'(\zeta',0, z')\equiv 0$,
$\Theta'(0,w',z')\equiv 0$. This implies in particular that the Segre
varieties ${\Cal S}_0 =\{(w,0)\: \v w \v < 2r \}$ and ${\Cal S}_0'
=\{(w',0)\: \v w'\v < 2r \}$ are straightened to the complex tangent
plane to $M$ at $0$ and that, if we develop $\bar\Theta$ and
$\bar\Theta'$ with respect to powers of $w$ and $w'$, then we can
write
$$
z=\xi+ i\sum_{\beta\in \N_*^{n-1}} 
{w}^\beta \ \bar\Theta_\beta(\zeta, \xi), \ \ \ \ \ 
z'=\xi'+ i\sum_{\beta\in \N_*^{n-1}} 
{w'}^\beta \ \bar\Theta_\beta'(\zeta', \xi').
\tag 2.4
$$
Here $\N_*^{n-1}:=\N^{n-1} \backslash \{0\}$, {\it i.e.} we mean that
the two above sums begin with a $w$ and $w'$ exponent of {\it
positive length}. It is now natural to set for notational convenience
$\bar\Theta_0(\zeta, \xi):=\xi$ and $\bar\Theta_0'(\zeta', \xi'):=\xi'$.
Although normal coordinates are in principle unnecessary, the
reduction to such normal coordinates will simplify a little the
presentation of all our formal calculations below.

\subhead 2.5.~Complexification of the map\endsubhead
Now, the map $h$ is by definition a $n$-vectorial formal power series
$h(t)=(h_1(t),\ldots,h_n(t))$ where $h_j(t)\in\C\dl t\dr$, $h_j(0)=0$
and $\hbox{det} \, \left( \partial h_j / \partial t_k(0)\right)_{1\leq
j,k\leq n}\neq 0$, which means that $h$ is formally invertible. 
This map yields
by extrinsic complexification a map $h^c=h^c(t,\tau)=(h(t),
\bar{h}(\tau))$ between the two complexification $({\Cal M},0)$ and $({\Cal
M'},0)$. In other words, if we denote $h=(g,f)\in
\C^{n-1}\times \C$ in accordance with the splitting of coordinates in
the target space, the assumption that $h^c({\Cal M})\subset_{\Cal F} 
{\Cal M}'$ reads as two equivalent fundamental equations\,:
$$
\Updownarrow
\left\{
\aligned
&
f(w,z)=\left[\bar{f}(\zeta,\xi)+i\bar{\Theta}'(g(w,z),\bar{g}(\zeta,\xi),
\bar{f}(\zeta,\xi))\right]_{\xi:=z-i\Theta(\zeta,w,z)},\\
&
\bar{f}(\zeta,\xi)=\left[f(w,z)-i\Theta'(\bar{g}(\zeta,\xi),g(w,z),
f(w,z))\right]_{z:=\xi+i\bar{\Theta}(w,\zeta,\xi)},
\endaligned\right.
\tag 2.6
$$
after replacing $\xi$ by $z-i\Theta(\zeta,w,z)$ in the first row and
$z$ by $\xi+i\bar{\Theta}(w,\zeta,\xi)$ in the second row. In fact,
these (equivalent) identities must be interpreted as {\it formal
identities} in the rings of {\it formal power series} $\C\dl
\zeta,w,z\dr$ and $ \C\dl w,\zeta,\xi\dr $ respectively. Of course,
according to~\thetag{2.3}, we can equally choose the coordinates
$(\zeta,w,z)$ or $(w,\zeta,\xi)$ over $\Cal M$. In symbolic notation,
we just write $h^c({\Cal M}, 0)
\subset_{\Cal F} ({\Cal M}',0)$ to mean the identities~\thetag{2.6}.

\subhead 2.7.~Conjugate equations, vector fields and the reflection 
function
\endsubhead
Let us also denote $r(t,\tau):=z-\xi-i\bar\Theta(w,\zeta, \xi)$,
$\bar r(\tau,t):=\xi-z+i\Theta(\zeta,w, z)$ and similarly
$r'(t',\tau'):= z'-\xi'-i\bar\Theta'(w',\zeta',
\xi')$, $\bar r'(\tau',t'):= \xi'-z'+i\Theta'(\zeta',w', z')$, so that
${\Cal M}=\{(t,\tau)\: r(t,\tau)=0\}$, ${\Cal M}'=\{(t',\tau')\:
r'(t',\tau')=0\}$ and the complexified Segre varieties are given by
${\Cal S}_{\tau_p}=\{ (t,\tau_p) \: r(t,\tau_p) =0\}\subset \Cal M$
for fixed $\tau_p$, and $\underline{\Cal S}_{t_p}=\{(t_p,\tau)\:
r(t_p, \tau)=0\}\subset \Cal M$ for fixed $t_p$ and similarly 
for ${\Cal S}_{\tau_{p'}'}'$, $\underline{\Cal S}_{t_{p'}'}'$ (again, 
the reader is referred to [9] for a complete exposition of the
geometry of complexified Segre varieties).
Finally, let us introduce the $(n-1)$ complexified (1,0) and (0,1) CR
vector fields tangent to $\Cal M$, that we will denote by ${\Cal L}
=({\Cal L}_1, \ldots, {\Cal L}_{n-1})$ and $\underline{\Cal L}=
(\underline{\Cal L}_1,\ldots,\underline{\Cal L}_{n-1})$, and which 
can be given in symbolic vectorial notation by
$$
{\Cal L}=\frac{\partial }{\partial w}+i\bar{\Theta}_w
(w,\zeta,\xi) \frac{\partial }{\partial z} 
\ \ \ \ \ \ \ \ {\text{\rm and}} 
\ \ \ \ \ \ \ \
\underline{\Cal L} =\frac{\partial }{\partial \zeta}-i\Theta_{\zeta}
(\zeta,w,z) \frac{\partial }{\partial \xi}.
\tag 2.8
$$
The {\it reflection function} $\Cal R_h' (t,\bar\nu')$, $t\in \C^n$,
$\bar\nu'=(\bar\lambda',\bar\mu')\in \C^{n-1} \times \C$, 
will be by definition the formal power series
$$
{\Cal R}_h'(t,\bar\nu')= {\Cal R}_h'(w,z,\bar\lambda', \bar\mu')
=\bar\mu'-f(w,z)+ i\sum_{\beta\in \N_*^{n-1}}
\bar{\lambda'}^\beta \ \Theta_\beta'(g(w,z), f(w,z)).
\tag 2.9
$$
We notice that this power series in fact 
belongs to the local ``hybrid'' ring
$\C\{\bar\nu'\}\dl t\dr$.

\head \S3. Minimality and holomorphic nondegeneracy \endhead

\subhead 3.1.~Two characterizations \endsubhead
At first, we need to remind the two explicit 
charac\-terizations of each one
of the main two assumptions of Theorem~1.2.
Let $M$ be a real analytic CR hypersurface given in {\it 
normal coordinates}
$(w,z)$ as above in eq.~\thetag{2.2}. 

\proclaim{Lemma~3.2}
\text{\rm ([3])}
The following properties are equivalent\,:
\roster
\item"{\bf (1)}" $\bar\Theta(w,\zeta,0)\not\equiv 0$.
\item"{\bf (2)}" ${\partial \bar\Theta \over \partial \zeta } 
(w,\zeta, 0)\not\equiv 0$.
\item"{\bf (3)}" $M$ is minimal at $0$.
\item"{\bf (4)}" The Segre variety $S_0$ is not contained in $M$.
\item"{\bf (5)}" The holomorphic map $\C^{2n-2} \ni (w,\zeta)\mapsto 
(w, i\bar\Theta (w,\zeta,0))\in \C^n$ has generic rank $n$.
\endroster
\endproclaim

\proclaim{Lemma~3.3}
\text{\rm{([2,3,16])}} If the coordinates $(w',z')$ are
normal as above, then the 
real analytic hypersurface $M'$ is holomorphically nondegenerate at
$0$ if and only if there exist $\beta^1, \ldots, \beta^{n-1}\in 
\N^{n-1}_*$, $\beta^n:=0$, such that
$$
\hbox{det} \ \left( 
{\partial \Theta_{\beta^i}'\over
\partial t_j'}(w',z')
\right)_{1\leq i,j\leq n}
\not\equiv 0 \ \ \ \ \ \hbox{in} \ \C\{w',z'\}.
\tag 3.4
$$
\endproclaim

\remark{Remark}
Since for $\beta=0$, we have $\Theta_0'(t')=z'$, we see
that~\thetag{3.4} holds if and only if $
\hbox{det} \ \left( {\partial \Theta_{\beta^i}'\over
\partial w_j'}(w',z') \right)_{1\leq i,j\leq n-1}
\not\equiv 0$. Further, we can precise the other
classical nondegeneracy conditions (I), (II) and (III) of \S1 (for
condition (IV), {\it see} Lemma~5.15)\,:

\proclaim{Lemma~3.5}
The following concrete characterizations hold in normal coordinates\,:
\roster
\item"{\bf (1)}"
$M'$ is Levi nondegenerate at $0$ if and only if the map $w'\mapsto
(\Theta_\beta(w',0))_{\v\beta\v=1}$ is immersive at $0$.
\item"{\bf (2)}"
$M'$ is finitely nondegenerate at $0$ if and only if there exists
$k_0\in \N_*$ such that the map $w'\mapsto
(\Theta_\beta(w',0))_{\v\beta\v\leq k_0}$ is immersive at $0$ for all
$k\geq k_0$. 
\item"{\bf (3)}"
$M'$ is essentially finite at $0$ if and only if there exists
$k_0\in \N_*$ such that the map $w'\mapsto
(\Theta_\beta(w',0))_{\v\beta\v\leq k_0}$ is finite at $0$ for all
$k\geq k_0$. 
\endroster
\endproclaim
\endremark

\subhead 3.6.~Switch of the assumptions \endsubhead
It is now easy to observe that the nondegeneracy conditions
upon $M$ transfer to $M'$ through $h$ and vice versa.

\proclaim{Lemma~3.7}
Let $h\: (M,0) \to_{\Cal F} (M',0)$ be a formal invertible CR
map between two real analytic hypersurfaces. Then
\roster
\item"{\bf 1)}" $(M,0)$ is minimal if and only if $(M',0)$ is minimal.
\item"{\bf 2)}" $(M,0)$ is holomorphically nondegenerate if and only if 
$(M',0)$ is holomorphically nondegenerate.
\endroster
\endproclaim

\demo{Proof}
We amit in the proof that minimality and holomorphic nondegeneracy are
biholomorphically invariant properties.  Let $N\in \N_*$ be
arbitrary. Since $h$ is invertible, after composing $h$ with a
biholomorphic and polynomial mapping $\Phi\: (M',0) \to (M'',0)$ which
cancels low order terms in the Taylor series $\Phi(h(t))=\sum_{\gamma\in
\N_*^{n-1}} h_\gamma t^\gamma$, we can achieve that $h(t)=t+ O(\vert
t\vert^N)$. Since the coordinates for $(M'',0)$ may be nonnormal, we
must compose $\Phi\circ h$ with a biholomorphism 
$\Psi\: (M'',0)\to (M''',0)$
which straightens the real analytic Levi-flat union of Segre varieties
$\bigcup_{\v x \v \leq r} S_{\overline{h(0,x)}}''$ into the real
hyperplane $\{y'''=0\}$ (this is how one constructs normal
coordinates). One can verify that $\Psi(t)=t+O(\v t\v^N)$ also.  Then
all terms of degre $\leq N$ in the power series of $\Theta'''$ coincide
with those of $\Theta$. Each one of the two characterizing properties
{\bf (1)} of Lemma~3.2 and~\thetag{3.4} of Lemma~3.3 is therefore
satisfied by $\Theta$ if and only if it is satisfied by $\Theta'''$.
\qed
\enddemo

\head \S4. Formal versus analytic \endhead

\subhead 4.1.~Approximation theorem\endsubhead
We collect here some useful statements from local analytic geometry
that we will repeatedly apply in the article. One of the essential arguments
in the proof of the main Theorem~2.1 rests on the existence of
analytic solutions arbitrarily close in the Krull topology to formal
solutions of some analytic equations, a fact which is known as {\it Artin's
approximation theorem}.  Let $\frak{m}(w)$ denote the maximal ideal of
the local ring $\C\dl w\dr$ of formal power series in $w\in \C^n$,
$n\in \N_*$. Here is the first of our three fundamental tools,
which will be used to get the Cauchy estimates
proving that the reflection function converges on the first Segre chain
({\it see}~Lemma~6.6).

\proclaim{Theorem~4.2} \text{\rm (Artin [1])}
Let $R(w,y)=0$, $R= (R_1,\ldots,R_J)$, where $w\in \C^n$, $y\in \C^m$,
$R_j\in\C\{w,y\}$, $R_j(0)=0$, be a converging system of \text{\rm
holomorphic} equations. Suppose
$\hat{g}(w)=(\hat{g}_1(w),\ldots,\hat{g}_m(w))$, $\hat{g}_k(w)\in
\C\dl w\dr$, $\hat{g}_k(0)=0$, are {\rm formal} 
power series which solve $R(w,\hat{g}(w))\equiv 0$ in $\C\dl
w\dr$. Then for every integer $N\in \N_*$, there exists a {\rm
convergent series solution} $g(w)= (g_1(w),\ldots,g_m(w))$, {\it i.e.}
satisfying $R(w,g(w))\equiv 0$, such that $g(w)\equiv \hat{g}(w) \
({\text{\rm mod}} \
\frak{m}(w)^N)$.
\endproclaim

\subhead 4.3.~Formal implies convergent\,: first recipe\endsubhead
The second tool will be used to prove that $h$ is
convergent on the second Segre chain, {\it i.e.} that
$h(w, i\bar\Theta(\zeta, w,0))\in \C\{w,\zeta\}$ ({\it see}
\S8).

\proclaim{Theorem 4.4}
Let $R(w,y)=0$, where $R=(R_1,\ldots,R_J)$, $w\in \C^n$, $y\in \C^m$,
$R_j\in\C\{w,y\}$, $R_j(0)=0$, be a system of \text{\rm
holomorphic} equations. Suppose that
$\hat{g}(w)=(\hat{g}_1(w),\ldots,\hat{g}_m(w))\in\C\dl w\dr^m$, 
$\hat{g}_k(0)=0$ are
\text{\rm formal} power series solving
$R(w,\hat{g}(w))\equiv 0$ in $\C\dl w\dr$. If $J\geq m$ and if there
exist $j_1,\ldots,j_m$ with $1\leq j_1 < j_2 < \cdots < j_m \leq J$ such
that
$$
{\text{\rm det}} \left(
\frac{\partial R_{j_k}}{\partial y_l}(w,\hat{g}(w))
\right)_{1\leq k,l \leq m} \not\equiv 0 
\ \ {\text{\rm in}} \ \ \C\dl w\dr,
\tag 4.5
$$
then the formal power series
$\hat{g}(w)\in \C\{w\}$ is in fact already convergent.
\endproclaim

\remark{Remark}
This theorem is a direct corollary of Artin's Theorem~4.2. The reader can
find an elementary proof of it for instance in \S12 of [10].
\endremark\medskip

\subhead 4.6.~Formal implies convergent\,: second recipe\endsubhead
The third statement will~be applied to the canonical map of the second
Segre chain, namely to the map $(w,\zeta)\mapsto (w, i\bar\Theta(\zeta,
w,0))$, which is of generic rank $n$ by Lemma~3.2 {\bf (5)}.

\proclaim{Theorem~4.7}
\text{\rm ([8])}
Let $a(y)\in \C\dl y\dr$, $y\in \C^\mu$, $a(0)=0$, be a formal power
series and assume that there exists a local holomorphic map $\varphi\:
(\C^\nu_x,0) \to (\C_y^\mu,0)$, of maximal generic rank $\mu$, {\it
i.e.} satisfying
$$
\exists \ j_1,\ldots,j_\mu, \ 1\leq j_1 < 
\cdots < j_\mu \leq \nu, \ \ s.t. \ \ \ 
\hbox{det} \ \left(
{\partial \varphi_k 
\over \partial x_{j_l}} 
(x)\right)_{1\leq k,l\leq
\mu} \not\equiv 0,
\tag 4.8
$$
and such that $a(\varphi(x)) \in \C\{x\}$ is convergent. Then 
$a(y)\in \C\{y\}$ is convergent.
\endproclaim

\subhead 4.9.~Application \endsubhead
We can now give an important application of Theorem 4.2\,: the 
Cauchy estimates for the convergence of the reflection function
come for free after one knows that all the formal power series
$\Theta_\beta'(h(w,z))\in \C\dl w,z\dr$ are convergent.

\proclaim{Lemma~4.10}
Assume that $h\:(M,0)\to_{\Cal F} (M',0)$ 
is a formal invertible CR mapping and that $M'$ is holomorphically
nondegenerate. Then the following properties are
equivalent\,:
\roster
\item"{\bf (1)}" 
$h(w,z)\in \C\{w,z\}^n$.
\item"{\bf (2)}" 
$\Cal R_h'(w,z,\bar\lambda,\bar\mu)\in 
\C\{w,z,\bar\lambda,\bar\mu\}$.
\item"{\bf (3)}" 
$\Theta_\beta'(h(w,z))\in \C\{w,z\}$, 
$\forall \ \beta \in \N^{n-1}$ 
and \ $\exists \ \varepsilon >0$ $\exists \ C >0$ such that
$\vert \Theta_\beta'(h(w,z))\vert \leq C^{\vert \beta\vert+1}$,
for all $(w,z)$ with 
$\v(w,z)\v < \varepsilon$ and all $\beta \in \N^{n-1}$.
\item"{\bf (4)}" 
$\Theta_\beta'(h(w,z))\in \C\{w,z\}$,
$\forall \ \beta \in \N^{n-1}$.
\endroster
\endproclaim

\demo{Proof} 
The implications {\bf (1)} $\Rightarrow$ {\bf (2)} $\Rightarrow$ {\bf
(3)} $\Rightarrow$ {\bf (4)} are straightforward.  On the other hand,
consider the implication {\bf (4)} $\Rightarrow$ {\bf (1)}.  By
assumption, there exist convergent power series
$\varphi_\beta'(w,z)\in \C\{w,z\}$ such that
$\Theta_\beta'(h(w,z))\equiv \varphi_\beta'(w,z)$ in 
$\C\dl w,z\dr$. It then follows that
$h(t)$ is convergent by an application of Theorem~4.4 with
$R_n(t,t'):=z'-\varphi_0(t)$ and $R_i(t,t'):=\Theta_{\beta^i}'(t')-
\varphi_{\beta^i}'(t)$, $1\leq i\leq n-1$ and where the multiindices 
$\beta^1,\ldots,\beta^{n-1}$ are chosen as in Lemma~3.3 (use the
property $\hbox{det}({\partial h_j\over
\partial t_k}(0))_{1\leq j,k\leq n}\neq 0$ and the composition
formula for Jacobian matrices to check~\thetag{4.5}).
\qed
\enddemo 

\demo{Proof of Proposition~1.5}
Let $\varphi'\: (t',u')\mapsto \exp (u'L')(t')=\varphi'(t',u')$ be the
local flow of the holomorphic vector field $L'=\sum_{k=1}^n a_k'(t')
\partial /\partial t_k'$ tangent to $M'$. Of course, this flow is
holomorphic with respect to $t'\in \C^n$ and $u'\in \C$, for $\v
t'\v$, $\vert u'\vert\leq \varepsilon$, $\varepsilon>0$. This flow
satisfies $\varphi'(t',0)\equiv t'$ and
$\partial_{u'}\varphi_k'(t',u')\equiv a_k'(\varphi'(t',u'))$. As
$L'\neq 0$, we have $\partial_{u'} \varphi'(t',u')\not\equiv 0$. We
can assume that $\partial_{u'}\varphi_1'(t',u')\not\equiv 0$. Let
$\varpi'(t')\in\C\dl t'\dr\backslash \C\{t'\}$, $\varpi'(0)=0$, be a
{\it nonconvergent} formal power series which satisfies further
$\partial_{u'} \varphi_1'(t',\varpi'(t'))\not\equiv 0$ in $\C\dl t'
\dr$ (there exist many of such). If the formal power series $h^\sharp\: 
t'\mapsto_{\Cal F}\varphi'(t',\varpi'(t'))$ would be convergent, then
$t'\mapsto_{\Cal F} \varpi'(t')$ would also be convergent, because of
Theorem~4.4, contrarily to the choice of $\varpi'$. Finally, $L'$
being tangent to $(M',0)$, it is clear that
$h^\sharp(M',0)\subset_{\Cal F} (M',0)$.
\qed
\enddemo

\head \S5. Classical reflection identities \endhead

\subhead 5.1.~The fundamental identities\endsubhead
In this paragraph, we start up the proof of our main Theorem~1.2 by
deriving the classical reflection identities. Thus let $\beta \in
\N_*^{n-1}$. By $\gamma \leq \beta$, we shall mean
$\gamma_1\leq\beta_1, \ldots, \gamma_{n-1}\leq \beta_{n-1}$. Denote
$\vert \beta\vert:=\beta_1+\cdots+\beta_{n-1}$ and $\underline{\Cal
L}^\beta:=\underline{\Cal L}_1^{\beta_1} \cdots
\underline{\Cal L}_{n-1}^{\beta_{n-1}}$. 
Then applying all these derivations of
any order ({\it i.e.} for each $\beta\in \N^{n-1}$)
to the identity $\bar r'(\bar h(\tau),h(t))$, {\it i.e.}  to
$$
\bar{f}(\zeta, \xi) \equiv f(w,z) -i\sumg \bar g(\zeta, \xi)^\gamma \
\Theta_\gamma'(g(w,z), f(w,z)),
\tag 5.2
$$
as $(w,z,\zeta,\xi)\in \Cal M$, it is well-known that we obtain an
infinite family of formal identities that we recollect here 
in an independent technical
statement (for the proof, {\it
see} [3,10]).

\proclaim{Lemma~5.3}
Let $h: (M,0) \to_{{\Cal F}} (M',0)$ be a formal invertible CR mapping
between ${\Cal C}^{\omega}$ hypersurfaces in $\C^n$. Then for every
$\beta\in \N^{n-1}_*$, there exists a collection of universal
polynomial $\underline{u}_{\beta,\gamma}$, $|\gamma|\leq |\beta|$ in
$(n-1)N_{n-1,|\beta|}$ variables, where $N_{k,l}:={(k+l)! \over k! \,
l!}$ and there exist holomorphic $\C$-valued functions
$\underline{\Omega}_{\beta}$ in $(2n-1+n N_{n,|\beta|})$ variables near
$0\times 0\times 0 \times
(\partial_{\xi}^{\alpha^1}\partial_{\zeta}^{\gamma^1}
\bar{h}(0))_{|\alpha^1|+ 
|\gamma^1| \leq |\beta|}$ in $\C^{n-1}\times \C^{n-1} \times \C \times
\C^{nN_{n,|\beta|}}$ such that the following identities
$$
\left\{
\aligned
&
{1\over \beta !} \partial_{\zeta'}^\beta \Theta' (
\bar g(\zeta, \xi), g(w,z),f(w,z)) =\\
&
= \Theta_\beta'(g(w,z), f(w,z))+\sumg 
{(\beta+\gamma)! \over \beta! \ \gamma!} \bar g(\zeta, \xi)^\gamma\
\Theta_{\beta+\gamma}'(g(w,z), f(w,z))\\
&
\equiv \sum_{\vert \gamma \vert \leq \vert \beta \vert} 
\frac{\underline{{\Cal L}}^{\gamma} \bar{f}(\zeta,\xi) \ \underline{u}_{\beta,\gamma}
((\underline{\Cal L}^{\delta} \bar{g}(\zeta,\xi))_{|\delta|\leq
|\beta|})}{\underline{\Delta}(w,\zeta,\xi)^{2|\beta|-1}}\\
&
=: \underline{\Omega}_{\beta}(w,\zeta,\xi,
(\partial_{\xi}^{\alpha^1}\partial_{\zeta}^{\gamma^1}\bar{h}
(\zeta,\xi))_{|\alpha^1|+|\gamma^1| \leq |\beta|})\\
&
=: \underline{\omega}_\beta(w,\zeta,\xi),
\endaligned\right.
\tag 5.4
$$
hold as formal power series in $\C\dl w,\zeta,\xi\dr$, where
$$
\left\{
\aligned
&\underline{\Delta}(w,\zeta,\xi)=\underline{\Delta}(w,z,\zeta,\xi)|_{z=
\xi+i\bar{\Theta}(w,\zeta,\xi)}:=
{\text{\rm det}}(\underline{\Cal L} \bar{g})=\\
&= {\text{\rm det}} \left(
\frac{\partial \bar{g}}{\partial \zeta}(\zeta,\xi)-
i\Theta_{\zeta}(\zeta,w,z)
\frac{\partial \bar{g}}{\partial \xi}(\zeta,\xi)\right)|_{z=\xi+
i\bar{\Theta}(w,\zeta,\xi)}.
\endaligned\right.
\tag 5.5
$$
\endproclaim
\remark{Remark}
The terms $\underline{\Omega}_{\beta}$, holomorphic in their
variables, arise after writing $\underline{\Cal
L}^{\delta}\bar{h}(\zeta,\xi)$ as $\chi_{\delta}(w,z,\zeta,\xi,
(\partial_{\xi}^{\alpha^1}\partial_{\zeta}^{
\gamma^1}\bar{h}(\zeta,\xi))_{|\alpha^1|+
|\gamma^1| \leq |\delta|})$ (by noticing that the coefficients of
$\underline{\Cal L}$ are analytic in $(w,z,\zeta,\xi)$) and by
replacing again $z$ by $\xi+i\bar{\Theta}(w,\zeta,\xi))$.
\endremark
 
\subhead 5.6.~Convergence over a uniform domain~\endsubhead
From this lemma which we have written down in the most explicit way,
we deduce the following useful observations. First, as we have by the
formal stabilization of Segre varieties $h(\{w=0\})\subset_{\Cal F}
\{w'=0\}$ and as $h$ is invertible, then it holds ${\text{\rm
det}}(\underline{\Cal L}\bar{g}(0))=\hbox{det}(
\partial g_j /\partial w_k(0))_{1\leq j,k\leq n-1}\neq 0$ also, whence
the rational term $1/\underline{\Delta}^{2|\beta|-1}\in \C\dl
w,\zeta,\xi\dr$ defines a true formal power series at the
origin. Putting now $(\zeta,
\xi)=(0,0)$ in eqs.~\thetag{5.5} and shrinking $r$ if necessary, we
then readily observe that $\underline{\Delta}^{1-2\vert
\beta \vert } (w,0,0) \in \Cal O ((r\Delta)^{n-1},\C)$, since
$\Theta_\zeta(0,w,0)
\in \Cal O ((r\Delta)^{n-1},\C)$ and since the terms 
$\partial_\zeta^{\gamma^1}
\bar g(0,0)$ for $\vert \gamma^1\vert =1$ and $\partial_\xi^1 \bar g
(0,0)$ are {\it constants}. Clearly, the numerator in the middle
identity \thetag{5.4} is also convergent in $(r\Delta)^{n-1}$ after
putting $(\zeta,\xi)=(0,0)$, and we deduce finally the following
important property\,:
$$
\underline{\Omega}_\beta(w,0,0,(\partial_\xi^{\alpha^1} 
\partial_\zeta^{\gamma^1} \bar h
(0,0))_{
\vert \alpha^1 \vert +\vert \gamma^1 \vert \leq \vert \beta \vert}
)\in \Cal O((r\Delta)^{n-1},\C),
\tag 5.7
$$
for all $\beta \in \N_*^{n-1}$. In other words, 
{\it the domains of convergence of
the $\underline{\omega}_\beta(w,0,0)$ are independent of $\beta$.}

\subhead 5.8.~Conjugate reflection identities\endsubhead
On the other hand, applying the same derivations 
$\underline{\Cal L}^\beta$'s to the conjugate identity 
$r'(h(t), \bar h(\tau))=0$, we would get another 
family of what we shall call {\it conjugate reflection identities}\,:
$$
0\equiv \underline{\Cal L}^\beta \bar f(\zeta, \xi)+i\sumg
g(w,z)^\gamma \ \underline{\Cal L}^\beta
(\bar\Theta_\gamma'(\bar g(\zeta, \xi), \bar f(\zeta, \xi))).
\tag 5.9
$$
But these equations furnish essentially no more information for the
reflection principle, because\,:

\proclaim{Lemma~5.10}
If $(t,\tau)\in \Cal M$, then
$$
\left<
\underline{\Cal L}^\beta(r'(h(t),\bar h(\tau)))=0, \ \forall
\beta \in\N^{n-1}\right> \ \Longleftrightarrow \ \left<
\underline{\Cal L}^\beta
(\bar r'(\bar h(\tau), h(t)))=0, \ \forall \beta \in \N^{n-1}
\right>.
\tag 5.11
$$
\endproclaim

\demo{Proof}
As the two equations for $\Cal M'$ are equivalent, there exists an
invertible formal series $\alpha(t,\tau)$ such that $r'(h(t), \bar
h(\tau))\equiv \alpha (t,\tau) \,
\bar r'(\bar h(\tau), h(t))$. Thus
$$
\left\{
\aligned
\underline{\Cal L}^\beta( r'(h(t), \bar h(\tau))) \ \equiv \
&
\alpha(t,\tau) \, \underline{\Cal L}^\beta (\bar r'(\bar h(\tau), h(t)))+\\
&
+\sum_{\gamma \leq \beta, \gamma\neq \beta} \alpha_\gamma^\beta
(t,\tau) \, \underline{\Cal L}^\gamma (\bar r'(\bar h(\tau), h(t))),
\endaligned\right.
\tag 5.12
$$
for some formal series $\alpha_\gamma^\beta(t,\tau)$ depending on the
derivatives of $\alpha(t,\tau)$. The implication ``$\Leftarrow$''
follows at once and the reverse implication is totally similar. $\square$
\enddemo

\subhead 5.13.~Heuristics \endsubhead
Nevertheless, in the last step of the proof of our main Theorem~1.2, the
equations~\thetag{5.9} above will be of crucial use, in place of the
equations~\thetag{5.4} which will happen to be unusable. The
explanation is the following. Whereas the jets
$(\partial_{\xi}^{\alpha^1}\partial_{\zeta}^{\gamma^1}\bar{h}(\zeta,
\xi))_{|\alpha^1|+ |\gamma^1| \leq |\beta|}$ of the mapping $\bar h$ 
cannot be seen directly to be convergent on the first Segre chain
$\Cal S_0^1:=\{(w,0)\}$, a convergence which would be a necessary fact
to be able to use formula \thetag{5.4} again in order to pass from the
first to the second Segre chain $\Cal S_2^0:=\{(w,i\bar
\Theta(w,\zeta,0))\}$ it will be possible -- fortunately! -- to show
in \S7 below that {\it the jets of the reflection function $\Cal R_h'$
itself converge on the first Segre chain}, namely that all the
derivatives $\underline{\Cal L}^\beta (\bar\Theta_\gamma'(\bar
g(\zeta, \xi), \bar f(\zeta,\xi)))$, restricted to the conjugate first
Segre chain $\underline{\Cal S}_0^1=\{(\zeta,0)\}$, converge. In
summary, we will only be able {\it a priori} to show that the jets of
$\Cal R_h'$ converge on the first Segre chain, and thus only the
equations~\thetag{5.9} will be usable in the next step, but not the
classical reflection identities~\thetag{5.4}.  {\it This shows
immediately why the conjugate reflection identities ~\thetag{5.9}
should be undertaken naturally in this context.}

\subhead 5.14.~The Segre-nondegenerate case \endsubhead
Nonetheless, in the Segre-nondegenerate case, which is less general
than the holomorphically nondegenerate case, we have been able to show
directly that the jets of $h$ converge on the first Segre chain ({\it
see} [10]), and so on by induction, without
using conjugate reflection identities. The explanation is simple\,:
in the Segre nondegenerate case, we have first the following
characterization, which shows that we can {\it separate} the $w'$
variables from the $z'$ variable\,:

\proclaim{Lemma~5.15}
The $\Cal C^\omega$ hypersurface $M'$, given in 
normal coordinates $(w',z')$, is Segre-nondegenerate at
$0$ if and only if there exist $\beta^1, \ldots, 
\beta^{n-1}\in \N_*^{n-1}$ such that
$$
\hbox{det} \ \left( 
{\partial \Theta_{\beta^i}'\over
\partial w_j'}(w',0)
\right)_{1\leq i,j\leq n-1}
\not\equiv 0 \ \ \ \ \ \hbox{in} \ \C\{w'\}.
\tag 5.16
$$
Also, $M'$ is holomorphically nondegenerate at $0$ if
it is Segre nondegenerate at $0$.
\endproclaim
 
\demo{Proof} 
In our normal coordinates, it follows that $\Cal S_{p'}=\Cal
S_0'=\{(w',0,0,0)\}$ and $\varphi_k'\v_{\Cal S_0'}\cong w'\mapsto
(\{\Theta_\beta'(w',0)\}_{\v \beta \v \leq k})$, whence the
rephrasing~\thetag{5.16} of~definition~(IV).  As we can take
$\beta^n=0$ in~\thetag{3.4} above, we see that the determinant
of~\thetag{3.4} does not vanish if~\thetag{5.16} holds. This proves
the promised implication $(IV)\Rightarrow (V)$.
\qed
\enddemo

\noindent
Thanks to this characterization, we can delineate an
analog to Lemma~4.10, whose proof goes exactly the same way\,:

\proclaim{Lemma~5.17}
Assume that $h$ is invertible, that $M$ is given in normal
coordinates 2.1 and that $M'$ is Segre nondegenerate. Then
the following properties are~equivalent
\roster
\item"{\bf (1)}" 
$h(w,0)\in \C\{w\}$.
\item"{\bf (2)}" 
$\Cal R_h'(w,0,\bar\lambda,\bar\mu)\in \C\{w,\bar\lambda,\bar\mu\}$.
\item"{\bf (3)}" 
$\Theta_\beta'(h(w,0))\in \C\{w\}$, $\forall \ \beta \in \N^{n-1}$ 
and \ $\exists \ \varepsilon >0$ $\exists \ C >0$ such that
$\vert \Theta_\beta'(h(w,0))\vert \leq C^{\vert \beta\vert+1}, \
\forall \ \vert w\vert < \varepsilon \ \ \forall \ \beta \in \N^{n-1}$.
\item"{\bf (4)}" 
$\Theta_\beta'(h(w,0))\in \C\{w\}$, $\forall \ \beta \in \N^{n-1}$.
\endroster
\endproclaim

\subhead 5.18.~Comment\endsubhead
In conclusion, in the Segre nondegenerate case (only) the convergence
of all the components $\Theta_\beta'(h)$ of the reflection mapping
{\it after restriction to the Segre variety $\Cal S_0=\{(w,0)\}$} is
{\it equivalent to the convergence of all the components of $h$}. The 
same property holds for jets. Thus, in the Segre nondegenerate case, 
one can use the classical reflection identities~\thetag{5.4}
(in which appear the jets of $\bar h$, see $\underline{\Omega}_\beta$)
by induction on the Segre chains [10]. This
is not so in the general holomorphically nondegenerate case, because
it can happen that~\thetag{3.4} holds whereas~\thetag{5.16} does not
hold, as shows the example~\thetag{1.10}. In substance,
one has therefore to use the {\it conjugate
reflection identities}. Now, the proof of our
main Theorem~1.2 will be subdivided in three steps, which will be achieved
in \S6, \S7 and \S8 below.

\head \S6. Convergence of the reflection function on $\Cal S_0^1$ \endhead

\subhead 6.1.~Examination of the reflection identities\endsubhead
The purpose of this paragraph is to prove as a first step that
the reflection function $\Cal R_h'$ converges on the first Segre chain 
$\Cal S_0=\{(w,0)\}$ or more precisely\,:

\proclaim{Lemma~6.2}
After perharps shrinking the radius 
$r>0$ of~\thetag{5.7}, the formal power series
$\Cal R_h'(w,0,\bar\lambda',\bar\mu')$ is holomorphic in
$(r\Delta)^{n-1}\times\{0\}\times (r\Delta)^n$.
\endproclaim

\demo{Proof}
We specify
the infinite family of identities~\thetag{5.2} (for $\beta=0$)
and~\thetag{5.4} (for $\beta\in \N_*^{n-1}$) on $\Cal S_0$, to obtain
first that $f(w,0)\equiv 0\in \C\{w\}$ and that for all $\beta\in
\N_*^{n-1}$
$$
\Theta_\beta'(g(w,0), f(w,0)) \equiv
\underline{\Omega}_\beta (w,0,0, 
(\partial_\xi^{\alpha^1}\partial_\zeta^{\gamma^1}\bar h(0,0))_{
\vert \alpha^1\vert + \vert \gamma^1 \vert \leq \vert \beta \vert })\ \ 
\in \ \C\{w\}.
\tag 6.3
$$
Furthermore, since by~\thetag{5.7} the $\underline{\Omega}_\beta$'s
converge for $\v w \v < r$ and $\zeta=\xi=0$, we have got
$\Theta_\beta'(g(w,0), f(w,0)) \in \Cal O((r\Delta)^{n-1}, \C)$,
$\forall \ \beta \in \N^{n-1}$. It remains to establish a Cauchy
estimate like in {\bf (3)} of Lemma~5.17. To this aim, we introduce
some notation.  We set $\varphi_0'(w,z):= f(w,z)$ and
$\varphi_\beta'(w,z):=
\Theta_\beta'(g(w,z), f(w,z))$ for all $\beta\in \N_*^{n-1}$.
By~\thetag{6.3}, we already know that all the series
$\varphi_{\beta}'(w,0)$ are holomorphic in $\{\v w\v < r\}$.  Thus, in
order to prove that the reflection function restricted to the first
Segre chain, namely that the series $=\Cal R_h'\vert_{\Cal S_0}
=\bar\mu'+i\sumb \bar{\lambda'}^\beta \ \varphi_\beta'(w,0)$ is 
convergent with respect to all its variables, we must establish a crucial
assertion.

\proclaim{Lemma~6.4}
After perharps shrinking $r>0$,
there exists a constant $C>0$ with
$$
\vert \varphi_\beta'(w,0)\vert
\leq C^{\vert \beta \vert +1}, \ \ \ \forall \ \v w\v < r,
\ \forall \ \beta \in \N^{n-1}.
\tag 6.5
$$
\endproclaim

\demo{Proof}
Actually, this Cauchy estimate will follow, by construction, from 
eq.~\thetag{6.3}
and from the property $\vert \Theta_\beta'(w',z')\vert \leq
{C'}^{\vert \beta \vert +1}$ when $(w',z')$ satisfy $\v
(w',z')\v < r'$ (the natural Cauchy estimate
for $\Theta'$), once we have proved the following
independent and important proposition, which is a rather
direct application Artin's
approximation Theorem 4.2\,:
\qed \enddemo

\proclaim{Lemma~6.6}
Let $w\in \C^\mu$, $\mu\in\N_*$, $\lambda(w)\in 
\C\dl w\dr^\nu$, $\lambda(0)=0$, $\nu\in \N_*$, and let
$\Xi_\beta(w,\lambda)\in \C\{w, \lambda\}$,
$\Xi_\beta(0,0)=0$, $\beta\in \N^m$, $m\in \N_*$,
be a collection of holomorphic functions satisfying
$$
\exists \ r >0 \ \ \exists \ C>0 \ \ s.t. \ \ \ 
\vert \Xi_\beta(w,\lambda)\vert \leq C^{\vert \beta \vert +1}, \ \ 
\forall \ \beta \in \N^m, \ \forall \ \v (w,\lambda)\v < r.
\tag 6.7
$$
Assume that $\Xi_\beta(w,\lambda(w))\in \Cal O((r \Delta)^\mu, \C)$,
$\forall \ \beta \in \N^m$ and put $\Phi_\beta(w):= 
\Xi_\beta(w,\lambda(w))$. Then
the following Cauchy inequalities are satisfied by the 
$\Phi_\beta$'s\,:
$$
\exists \ 0 < r_1 \leq r,
\ \exists \ C_1 >0 \ \ s.t. \ \ \vert \Phi_\beta (w)\vert
\leq C_1^{\vert \beta \vert +1}, \ \forall \ \beta \in
\N^m, \ \forall \ \v w\v < r_1.
\tag 6.8
$$
\endproclaim

\demo{Proof}
We set $R_\beta(w,\lambda):= \Xi_\beta(w,\lambda)-\Phi_\beta(w)$. Then
$R_\beta\in {\Cal O}(\{\v (w,\lambda)\v < r \}, \C)$.
By noetherianity, we can assume that a
finite subfamily $(R_\beta)_{\vert \beta \vert \leq \kappa_0}$
generates the ideal $(R_\beta)_{\beta \in \N^m}$, for some
$\kappa_0\in \N_*$ large enough.
Applying now Theorem~4.2 to the collection of equations
$R_\beta(w,\lambda)=0$, $\vert \beta \vert \leq \kappa_0$, of which a
formal solution $\lambda(w)$ exists by assumption, we get that there
exists a convergent solution $\lambda_1(w)\in \C\{w\}^\nu$ vanishing
at the origin, {\it i.e.} some $\lambda_1(w)\in \Cal O((r_1
\Delta)^\mu, \C^\nu)$, for some $0< r_1 \leq r$, with
$\lambda_1(0)=0$, which satisfies $R_\beta(w,
\lambda_1(w))\equiv 0$, $\forall \ \vert \beta \vert \leq
\kappa_0$. This implies that $R_\beta(w, \lambda_1(w))\equiv 0$,
$\forall \ \beta \in \N^m$. Now, we have obtained
$$
\Xi_\beta(w,\lambda(w)) \equiv \Phi_\beta(w)\equiv 
\Xi_\beta(w,\lambda_1(w)), 
\ \ \ \ \ \forall \ \beta \in \N^m.
\tag 6.9
$$
The composition formula for analytic function then yields at once $\vert
\Xi_\beta (w,\lambda_1(w))\vert \leq C_1^{\vert \beta \vert +1}$ for
$\vert w \vert < r_1$, after perharps shrinking once more this
positive number $r_1$ in order that $\vert \lambda_1(w)
\vert < r/2$ if $\vert w \vert < r_1$. Thanks to
eq. (6.9), this gives the desired inequality for $\Phi_\beta(w)$.
The Proof of Lemmas~6.7 and~6.2 are thus complete now.
\qed
\enddemo
\enddemo

\head \S7. Convergence of the jets of the reflection function on 
$\Cal S_0^1$ \endhead

\subhead 7.1.~Transversal differentiation of the reflection 
identities ~\endsubhead
The next step in our proof consists in showing that all the 
jets of the reflection function converge on the first Segre chain
$\Cal S_0$, or more precisely\,:

\proclaim{Lemma~7.2}
For all $\alpha\in \N$ and all $\gamma \in \N^{n-1}$, we have
$$
\left[ \partial_z^\alpha \partial_w^\gamma \Cal R_h'(w,z,
\bar\lambda ,\bar\mu)\right]
\vert_{z=0} \ \ \in \ \C\{w,\bar\lambda,\bar\mu\}.
\tag 7.3
$$
Equivalently, $\forall \ \alpha \in \N$,
$\forall \ \gamma \in \N^{n-1}$, $\exists \
r(\alpha,\gamma)>0$, $\exists \ C(\alpha,\gamma) >0$ 
such that
$$
\vert [\partial_z^\alpha\partial_w^\gamma \varphi_\beta'(w,z)]
\vert_{z=0} \vert 
\leq C(\alpha, \gamma)^{\vert \beta \vert +1} \ \ if \ \ \vert w \vert
< r(\alpha,\gamma), \ \ \forall \ \beta \in \N_*^{n-1}.
\tag 7.4
$$
\endproclaim

\remark{Remark}
Fortunately, the fact that $r(\alpha,\gamma)$ depends on $\alpha$ and
$\gamma$ will cause no particular obstruction for the achievement of
the last third step in \S8 below. We believe however that this
dependence should be avoided, but we get no immediate control of
$r(\alpha,\gamma)$ as $\alpha+\vert \gamma \vert \to
\infty$, in our proof -- although it can be seen by induction that $
[\partial_z^\alpha\partial_w^\gamma \varphi_\beta'(w,z)]\vert_{z=0}\in
\Cal O((r\Delta)^{n-1}, \C)$ ({\it cf.} 
the proof of Lemma~7.2 below).
\endremark

\demo{Proof}
If we denote by ${\Cal E}_{\alpha, \gamma}$ the statement of the
lemma, then it is clear that
$$
\Cal E_{\alpha, 0} \ \Rightarrow \ (\Cal E_{\alpha, \gamma} 
\ \ \forall \ \gamma \in 
\N^{n-1}).
\tag 7.5
$$
It suffices therefore to establish the truth of 
$\Cal E_{\alpha, 0}$ for all $\alpha \in \N$.
Let us first establish that 
$\partial_z^\alpha\vert_{z=0} [\varphi_\beta'(w,z)]\in 
\Cal O((r \Delta)^{n-1}, \C)$, $\forall \ \alpha \in \N$,
$\forall \ \beta \in \N^{n-1}$.
To this aim, we specify the variables $(w,z,\zeta,\xi):=
(w,z,0,z)\in \Cal M$
(because $\Theta(0,w,z) \equiv 0$) in the equations~\thetag{5.2}
and~\thetag{5.4} to obtain firstly
$$
\bar f(0, z) \equiv f(w,z) - i\sumg \bar g(0, z)^\gamma \ 
\Theta_\gamma'(g(w,z), f(w,z))
\tag 7.6
$$
and secondly the following infinite number of relations\,:
$$
\left\{
\aligned
&
\underline{\Omega}_\beta(w,0,z,(
\partial_\xi^{\alpha^1} \partial_{\zeta}^{\gamma^1} 
\bar h(0,z))_{\vert \alpha^1 \vert +\vert \gamma^1 
\vert \leq \vert \beta \vert})
\equiv\\
&
\equiv \Theta_\beta'(g(w,z), f(w,z))+\sumg 
{(\beta+\gamma)!\over \beta!\ \gamma!} \
\bar g(0,z)^\gamma \ \Theta_{\beta+\gamma}'(g(w,z), f(w,z)).
\endaligned\right.
\tag 7.7
$$
Essentially, the game will consist in differentiating the
equalities~\thetag{7.6} and~\thetag{7.7} with respect to $z$ at $0$ up
to arbitrary order $\alpha$, in the aim to obtain new identities which
will yield $\partial^\alpha_z\vert_{z=0} [ \Theta_\beta'(g(w,z),
f(w,z))]\in \Cal O((r \Delta)^{n-1}, \C)$, $\forall \ 
\beta \in \N^{n-1}$, $\forall \ \alpha \in \N$, by an induction 
process of ``trigonal type''. Let us complete this informal description.
To begin with, for $\alpha
= 1$, after applying the derivation operator $\partial_z^1\vert_{z=0}$
to eqs.~\thetag{7.6} and~\thetag{7.7}, we get immediately
$$
\partial_z^1 f(w,0) \equiv \partial_z^1\bar f(0,0)+ i 
\sum_{j=1}^{n-1} [\partial \bar g_j/\partial z] (0,0) \ 
\Theta_\gamma'(g(w,0), f(w,0))\in \C\{w\},
\tag 7.8
$$
since $\bar g(0,0)=0$ (so $\partial_z^1[\bar g(0,z)^\gamma]_{z=0}=0$ for
$\vert \gamma \vert \geq 2$), and
$$
\left\{
\aligned
\partial_z^1 \Theta_\beta'(g(w,0), f(w,0))
& 
= [\partial_z^1[\underline{\Omega}_\beta
(w,0,z, (\partial_\xi^{\alpha^1}\partial_\zeta^{\gamma^1} \bar h(0,z))_{
\vert \alpha^1 \vert+\vert \gamma^1 \vert 
\leq \vert \beta \vert})]]\vert_{z=0}-\\
&
-\sum_{\vert \gamma \vert =1} {(\beta+\gamma)!\over \beta!\ \gamma!} 
\partial_z^1 \bar g(0,0)^\gamma \ \Theta_{\beta+\gamma}'(g(w,0), f(w,0)),
\endaligned\right.
\tag 7.9
$$
making the slight abuse of notation $\partial_z^1 \chi(w,0)$ instead
of writing $\partial_z^1 \vert_{z=0} [\chi(w,z)]$ for any formal power
series $\chi(w,z)\in \C\dl w,z\dr$. For instance, $\partial_z^1
\bar g(0,0)^\gamma$ significates 
$[\partial_z^1 (\bar g(0,z)^\gamma)]\vert_{z=0}=
\sum_{k=1}^{n-1} \gamma_k \partial_z^1 \bar g_k(0,0) [
\bar g(0,0)^{\gamma_1}\cdots
\bar g_k(0,0)^{\gamma_k-1} \cdots \bar g_{n-1}(0,0)^{\gamma_{n-1}}]$.
All these expressions are convergent, because we know already (thanks
to the first step) that $\Theta_\beta'(g(w,0), f(w,0))\in \C\{w\}$
$\forall \ \beta \in \N_*^{n-1}$ (and even $\Theta_\beta'(g(w,0),
f(w,0))\in \Cal O((r \Delta)^{n-1}, \C)$) and because the derivative
$\partial_z^1 \vert_{z=0}(\underline{\Omega}_\beta)$ can be expressed
(thanks to the chain rule) in terms of the derivatives $\partial
\underline{\Omega}_\beta / \partial z$, 
in terms of the derivatives $\partial
\underline{\Omega}_\beta / (\partial^{\alpha^1} \partial^{\gamma^1})$
(considering $ \partial^{\alpha^1}
\partial^{\gamma^1}$ as independent variables), 
and in terms of the derivatives $\partial_z^1\vert_{z=0}
(\partial_\xi^{\alpha^1}\partial_\zeta^{\gamma^1}
\bar h(0,z))$, all taken at $z=0$, which are terms 
obviously converging and even which belong to the space $\Cal
O((r\Delta)^{n-1}, \C)$.
Thus, we have got that $\partial_z^1 \varphi_\beta'(w,0)\in 
\Cal O((r \Delta)^{n-1}, \C)$, $\forall \ \beta \in \N^{n-1}$
(including $\beta=0$).
More generally, for arbitrary $\alpha\in \N$ and $\beta\in \N_*^{n-1}$, 
we observe readily that $\partial_z^\alpha\vert_{z=0} [\bar f(0, z)]$ is
constant and, for the same reasons as explained above, that
$$
\partial_z^\alpha\vert_{z=0} \left[
\underline{\Omega}_\beta (w,0,z, \right.
\left. (\partial_\xi^{\alpha^1}\partial\zeta^{\gamma^1} 
\bar h(0,z))_{
\vert \alpha^1 \vert+\vert \gamma^1 \vert 
\leq \vert \beta \vert})\right]\in
\Cal O((r \Delta)^{n-1}, \C).
\tag 7.10
$$
We can use this observation to perform a ``trigonal'' induction as
follows.  Let $\alpha_0\in \N_*$ and suppose by induction that
$\partial^\alpha_z
\varphi_\beta'(w,0) \in \Cal O((r \Delta)^{n-1}, \C)$ 
$\forall \ \alpha \leq \alpha_0$, $\forall \ \beta \in \N_*^{n-1}$.
Then applying the derivation $\partial^{\alpha_0+1}_z\vert_{z=0}$
to~\thetag{7.6}, developing the expression according to Leibniz'
formula and using the fact that
$\partial^{\alpha_0+1}_z\vert_{z=0} [\bar g(0,z)^\gamma]=0$ for all
$\vert \gamma \vert \geq \alpha_0+2$, we get the expression\,:
$$
\left\{
\aligned
&
\partial_z^{\alpha_0+1}f(w,0) 
\equiv \partial_z^{\alpha_0+1}
\bar f(0,0) +i\sum_{0<\vert \gamma\vert \leq \alpha_0+1}\\
&
\sum_{\kappa=1}^{\alpha_0+1} 
{(\alpha_0+1)!\over \kappa! (\alpha_0+1-\kappa)!} \
\partial_z^\kappa \bar g(0,0)^\gamma \ \ \partial_z^{\alpha_0+1-\kappa}
\varphi_\gamma'(w,0).
\endaligned\right.
\tag 7.11
$$
Now, thanks to the induction assumption and {\it because the order of
derivation in the expression $\partial_z^{\alpha_0+1-\kappa}
\varphi_\gamma'(w,0)$ for $1\leq \kappa \leq \kappa_0+1$
is less or equal to $\alpha_0$}, we obtain that
this expression~\thetag{7.11} belongs to 
$\Cal O((r \Delta)^{n-1}, \C)$. Concerning
the differentiation of~\thetag{7.9} with
respect to $z$, we also get that the term
$$
\left\{
\aligned
&
\partial_z^{\alpha_0+1} \varphi_\beta'(w,0) 
\equiv \partial_z^{\alpha_0+1}
\underline{\omega}_\beta(w,0,0)-\sum_{0< \vert \gamma \vert \leq \alpha_0+1} 
{(\beta+\gamma)!\over \beta!\ \gamma!} \\
&
\left(
\sum_{\kappa=1}^{\alpha_0+1} 
{(\alpha_0+1)! \over \kappa ! (\alpha_0+1 -\kappa)!} \
\partial_z^\kappa \bar g(0,0)^\gamma \ \ 
\partial_z^{\alpha_0+1-\kappa} \varphi_{\beta+\gamma}'(w,0)\right)
\endaligned\right.
\tag 7.12
$$
belongs to $\Cal O((r \Delta)^{n-1}, \C)$. Again, the important fact
is that in the sum $\sum_{\kappa=1}^{
\alpha_0+1}$, only the derivations 
$\partial_z^\alpha \varphi_\beta'(w,0)$ 
for $0\leq \alpha \leq \alpha_0$ occur.
In summary, we have shown that $\partial_z\varphi_\beta(w,0)$
is convergent for all $\alpha\in\N$.
\enddemo

\subhead 7.13.~Intermezzo \endsubhead
The induction process can be said to be of ``trigonal type'' because
we are dealing with the infinite collection of identities~\thetag{5.4}
which can be interpreted as a linear system $Y=AX$, where $X$ denotes
the unknown $(\Theta_\beta')_{\beta\in\N^{n-1}}$ and $A$ is an
infinite trigonal matrix, as shows an examination
of~\thetag{5.4}. Further, when we consider the jets, we still get a
trigonal system. {\it The main point is that after restriction to the
first Segre chain $\{\zeta=\xi=0\}$, this trigonal system becomes
diagonal (or with only finitely many nonzero elements after applying
$\partial_z^\alpha$), but this crucial simplifying property fails to
be satisfied after passing to the next Segre chains.} To be honest, we
should recognize that the proof we are conducting here unfortunately
fails (for this reason) to be generalizable to higher
codimension$\ldots$ However, an important natural idea will appear
during the course of the proof, namely the appearance of the natural
(and new) {\it conjugate reflection identities}~\thetag{5.9} which we
will heavily use in \S8 below. For reasons of symmetry, we have
naturally wondered whether they can be exploited more deeply. A
complete investigation is contained in our subsequent work on the
subject (quoted in \S1.12).

\subhead 7.14.~End of proof of Lemma~7.2\endsubhead
It remains to show that there exist constants $r(\alpha)>0$,
$C(\alpha)>0$ such that the estimate~\thetag{7.4} holds for $(\alpha,
\gamma)=(\alpha,0)$\,: $\vert \partial_z^\alpha \varphi_\beta'(w,0)
\vert \leq C(\alpha)^{\vert \beta \vert +1}$ if $\vert w \vert <
r(\alpha)$, $\forall \ \beta\in \N_*^{n-1}$.  To this aim, we shall
apply Lemma~6.6 with the suitable functions and variables.  First, it
is clear that there exist universal polynomials
\footnote{
The explicit formula in dimension one 
for the derivative of a composition
${d^n\over dx^n} (\psi\circ \phi (x))=
(\psi\circ\phi)^{(n)}(x)$
is known as {\it Faa di Bruno's formula},
(one of the favorite students of Cauchy):
$$
\aligned
{1\over n!} (\psi\circ \phi)^{(n)}(x)
=
&
\sum_{\alpha_1+2\alpha_2+\cdots
+n\alpha_n=n}
{1\over
\alpha_1! \, \alpha_2! \cdots \alpha_n! \,
(1!)^{\alpha_1} \, (2!)^{\alpha_2} \cdots (n!)^{\alpha_n}}\\
\times 
&
(\phi'(x))^{\alpha_1} \, (\phi''(x))^{\alpha_2} \cdots 
(\phi^{(n)}(x))^{\alpha_n} \ 
\psi^{(\alpha_1+\alpha_2+\cdots+\alpha_n)}(\phi(x)).
\endaligned
$$
(We ignore wether similar formulas in several variable 
are known or attributed.)
} 
such that the following composite derivatives can be written
$$
\partial_z^\alpha [\Theta_\beta'(h(w,z))]=
P_\alpha(\nabla^{*\alpha}_z h(w,z)), (
\nabla_{t'}^{*\alpha} \Theta_\beta')(h(w,z))),
\tag 7.15
$$
where the $n\alpha$-tuple 
$\nabla^{*\alpha} h(w,z):=(
((\partial_z^k h_1(w,z),\ldots,\partial_z^kh_n(w,z))_{
1\leq k\leq \alpha})$ and
the $({(\alpha+n)!\over \alpha! \ n!}-1)$-tuple 
$\nabla^{*\alpha}_{t'}\Theta_\beta'(t'):=
((\partial_{t'}^\beta \Theta_\beta'(t'))_{
1\leq \vert \beta \vert\leq \alpha}$.
We now consider these polynomials as holomorphic functions
$G_\beta^\alpha=G_\beta^\alpha(\nabla_z^\alpha h)$ of the
$n(\alpha+1)$ variables $\nabla^\alpha_z h=
((\partial^k_z h_1, \ldots, \partial^k_zh_n)_{0
\leq k \leq \alpha})$ which satisfy, by 
eq.~\thetag{7.15}\,:
$$
\partial_z^\alpha \Theta_\beta'(h(w,z))= 
G_\beta^\alpha(\nabla_z^\alpha h(w,z))
=G_\beta^\alpha (\nabla_z^\alpha h)\vert_{
\nabla^\alpha_z h:= \nabla_z^\alpha h(w,z)},
\tag 7.16
$$
where the $n\alpha$-tuple $\nabla_z^\alpha h(w,z)=
(\partial_z^k h_j(w,z))_{0\leq k\leq \alpha}^{1\leq j\leq n}$ and
$\nabla_z^\alpha h:=(\partial_z^k h_j)_{0\leq k\leq 
\alpha}^{1\leq j\leq n}$
are $n(\alpha+1)$ {\it independent variables} as we have just said above.
Obviously, these functions $G_\beta^\alpha(\nabla_z^\alpha h)$ 
satisfy an estimate of the form
$\vert G_\beta^\alpha(\nabla_z^\alpha h)\vert \leq C(\alpha)^{
\vert \beta \vert +1}$ if $\vert \nabla_z^\alpha h \vert < r$, 
because the functions $\nabla_{t'}^{*\alpha} \Theta_\beta'(t')$ 
satisfy an estimate of the form
$\vert \nabla_{t'}^{*\alpha} \Theta_\beta'(t')\vert\leq
C'(\alpha)^{\vert \beta \vert +1}$ if
$\vert t' \vert < r'$, for some constants $C'(\alpha)>0$, 
$r'>0$, and because we have 
$P_\alpha( \nabla_z^\alpha h,0)\equiv 0$.
We already know that there exist holomorphic functions 
$\chi_\beta^\alpha(w)=\partial_z^\alpha \varphi_\beta'(w,0)$
in $\{\vert w\vert < r\}$ indexed by $\beta\in \N_*^{n-1}$ 
such that the following formal identity holds:
$$
G_\beta^\alpha(\nabla_z^\alpha h(w,0))= 
\partial_z^\alpha \varphi_\beta'(w,0) =
\chi_\beta^\alpha (w) \ \ \ \hbox{in}
\ \C\dl w \dr.
\tag 7.17
$$
Now, a direct application of Lemma~6.6 yields the desired estimate\,:
$$
\vert \partial_z^\alpha \varphi_\beta'(w,0)\vert \leq
C(\alpha)^{\vert \beta \vert +1}\ \ \hbox{if} \ \ 
\vert w \vert < r(\alpha).
\tag 7.15
$$
Thus, we have completed the proof of Lemma~7.2.
\qed

\remark{Important remark}
When $\alpha\to \infty$, the number $(n+1)\alpha$ of variables in 
$\nabla_z^\alpha h$ also becomes infinite. Thus, at each step we apply
Artin's Theorem in Lemma~6.6, the $r(\alpha)$ may shrink
and go to zero as $\alpha \to \infty$.
\endremark 

\head \S8. Convergence of the mapping \endhead

\subhead 8.1.~Jump to the second Segre chain\endsubhead
We now complete the final third step by establishing that the power
series $h(t)$ is convergent in a neighborhood of $0$.  Let
$\underline{\Cal S}_0^2=\{\exp w \Cal L (\exp \zeta \underline{\Cal L}
(0)) \:\vert w\vert < r, \vert \zeta \vert < r\}$ be the second conjugate
Segre chain [9], or equivalently in our normal coordinates
$\underline{\Cal S}_0^2=\{(w, i\bar\Theta(w,\zeta,0), \zeta, 0) \:
\vert w \vert < r, \vert \zeta \vert < r\}$. We shall prove that
the map $h^c$ is convergent on $\underline{\Cal S}_0^2$.
More precisely\,:

\proclaim{Lemma~8.2}
The formal power series
$h(w, i\bar\Theta(\zeta, w,0))\in \C\{w,\zeta\}^n$ is convergent.
\endproclaim

\noindent
From this lemma, we see now how to achieve the proof of our Theorem~1.2\,:

\proclaim{Corollary~8.3}
Then the formal power series $h(w,z)\in \C\{w,z\}^n$ is
convergent.
\endproclaim

\demo{Proof}
We just apply Theorem~4.7, taking into account {\bf (5)} of
Lemma~3.2. $\square$
\enddemo

\demo{Proof of Lemma~8.2}
Thus, we have to show that $h(w,i\bar{\Theta}(\zeta, w,0))\in
\C\{w,\zeta\}$.  To this aim, we consider the conjugate reflection
identities~\thetag{5.2} and~\thetag{5.9} for various $\beta\in
\N_*^{n-1}$ after specifying them over 
$\underline{\Cal S}_0^2$, {\it i.e.}  after
setting $(w,z,\zeta,\xi):= (w,i\bar\Theta(\zeta, w,0), \zeta, 0)\in
\Cal M$, which we may write explicitely as follows
$$
\left\{
\aligned
&
\bar f(\zeta, 0) \equiv f(w, i\bar\Theta(w,\zeta,0))-i\sumg
\bar g(\zeta, 0)^\gamma \ \Theta_\gamma'(h(w,i\bar\Theta (w,\zeta,0)),\\
&
0\equiv [\underline{\Cal L}^\beta \bar f(\zeta, \xi)]_{\xi=0} +i\sumg
g(w,i\bar\Theta(w,\zeta,0))^\gamma \
[\underline{\Cal L}^\beta (\bar\Theta_\gamma'(\bar h(\zeta,\xi)))]_{\xi=0},
\endaligned\right.
\tag 8.4
$$
for all $\beta\in \N_*^{n-1}$.
Let now $\kappa_0\in \N_*$ be an integer larger than the supremum of
the lengths of some multiindices $\beta^i$'s, $1\leq i\leq n-1$,
satisfying the determinant property stated in eq.~\thetag{3.4} of
Lemma~3.3, {\it i.e.} $\kappa_0 \geq \sup_{1\leq i\leq n-1}
\vert \beta^i\vert$. According to Lemma~7.2, if we consider
the equations~\thetag{8.4} only for a finite number of $\beta$'s, say for
$\vert \beta \vert \leq \kappa_0$, there will exist a positive
number $r_1>0$ with
$r_1 \leq r$ and a constant $C_1>0$ such that each power series
$[\underline{\Cal L}^\beta (\bar\Theta_\gamma '(\bar h(\zeta,
\xi)))]_{\xi=0}=:
\chi_\gamma^\beta(w,\zeta)$ is holomorphic in the polydisc
$\{\vert w\vert, \vert \zeta \vert < r_1\}$ and satisfies the Cauchy
estimate $\vert
\chi_\gamma^\beta(w,\zeta)\vert \leq C_1^{\vert \gamma \vert +1}$ when
$\vert (w,\zeta)\vert < r_1$.  We can now represent the
eqs.~\thetag{8.4} under the brief form
$$
s_\beta(w,\zeta, h(w,i\bar{\Theta}(w,\zeta,0)))\equiv 0, \ \ \ 
\vert \beta \vert\leq \kappa_0,
\tag 8.5
$$
where the holomorphic functions $s_\beta=s_\beta(w,\zeta,t')$ are
simply defined by replacing the terms $[\underline{\Cal L}^\beta
(\bar\Theta_\gamma'(\bar h(\zeta, \xi)))]_{\xi=0}$ by
$\chi_\gamma^\beta(w,\zeta)$ in eqs.~\thetag{8.4}, so that the power
series $s_\beta$ converge in the set $\{\vert w\vert, \vert \zeta
\vert, \vert t' \vert < r_1\}$.  The goal is now to apply 
Theorem~4.4 to the collection of equations~\thetag{8.5} in
order to deduce that $h(w, i\bar\Theta(w,\zeta,0))\in \C\{w,\zeta\}$.

\remark{Remark}
As noted in the introduction, another (more powerful) idea would
be to apply the Artin Approximation Theorem~4.2 to the 
equations~\thetag{8.5} to deduce the existence of a {\it converging}
solution $H(w,\zeta)$ and then to deduce that the reflection function
itself converges on the second Segre chain (which is in fact quite 
easy using Lemma~5.10). It's a pitty that in September-October 1999, 
while writting this \S8, and being highly conscient that
the conjugate reflection identities~\thetag{5.9} should be capital,
I missed the simplest idea.
\endremark

\smallskip

First, we make a precise choice of the $\beta^i\in \N_*^{n-1}$ arising
in Lemma~3.3. We set $\underline{\beta}^n=0$ and, for
$1\leq i\leq n-1$, let $\underline{\beta}^i$ be the 
infimum of all the multiindices $\beta\in \N_*^{n-1}$ satisfying
$\beta> \underline{\beta}^{i+1} > \cdots > \underline{\beta}^{n}$ for the
natural lexicographic order on $\N^{n-1}$, 
and such that an $(n-i+1)\times (n-i+1)$ 
minor of the $n\times (n-i+1)$ matrix
$$
{\Cal M\Cal A\Cal T}_{\beta, \underline\beta^{i+1},\ldots, 
\underline{\beta}^n}(t'):=
\left(
{\partial \Theta_\beta'\over \partial t_j'}(t') \ \
{\partial \Theta_{\underline{\beta}^{i+1}}'\over \partial t_j'}(t') \ \ 
\cdots \ \
{\partial \Theta_{\underline{\beta}^{n}}'\over \partial t_j'}(t')
\right)_{1\leq j\leq n}
\tag 8.6
$$
does not vanish identically as a holomorphic function of $t'\in
\C^n$. We thus have $\hbox{det} \left( {\partial
\Theta_{\underline{\beta}^i}' \over
\partial t_j'}(t') \right)_{1\leq i,j\leq n}\not\equiv 0$ in
$\C\{t'\}$. Concerning the choice of $\kappa_0$, we also require that
$$
\kappa_0 \geq \inf \{k\in \N \: 
\hbox{det} \left(
{\partial \Theta_{\underline{\beta}^i}'\over \partial
t_j'} (h(w,i\bar\Theta(w,\zeta,0))) 
\right)_{1\leq i,j\leq n}\not\in \ 
\frak{m} (\zeta)^k \ \C\dl w,\zeta \dr\},
\tag 8.7
$$
where $\frak{m}(\zeta)$ is the maximal ideal of $\C\dl \zeta \dr$.  We
can choose such a finite $\kappa_0$, because we know already that the
determinant in
\thetag{8.7} does not vanish identically (this fact can be easily 
checked, after looking at the composition formula for Jacobians,
because, in view of Lemma~3.3, the determinant ~\thetag{8.6} for
$(\underline{\beta}^1,\ldots,\underline{\beta}^n)$ does not vanish
identically and because the determinant $\hbox{det} \left({\partial
h_j\over
\partial t_k}(w,i\bar\Theta(w,\zeta,0))\right)_{1\leq j,k\leq n}$
does not vanish identically in view of the invertibility assumption on
$h$ and in view of the minimality criterion Lemma~3.2 {\bf (5)}).
Thus, after these choices are made, in order to finish the proof by an
application of Theorem~4.4, it will suffice to show that\,:

\proclaim{Lemma~8.8}
There exist $\beta^1, \ldots, \beta^{n-1}, \beta^n(=0)\in \N^{n-1}$
with $\vert \beta^i\vert \leq 2\kappa_0$ such that 
$$
\hbox{det} \left( {\partial s_{\beta^i}\over \partial t_j'}
(w,\zeta, h(w,i\bar\Theta(w,\zeta,0)))\right)_{1\leq i,j\leq n}
\not\equiv 0 \ \ \hbox{in} \ \C\dl w,\zeta\dr.
\tag 8.9
$$
\endproclaim

\demo{Proof} To this aim, we introduce some new
power series. We set\,:
$$
R_\beta(w,z,\zeta,t'):= \underline{\Cal L}^\beta \bar f(\zeta, \xi)+
i\sumg \underline{\Cal L}^\beta (\bar g(\zeta, \xi)^\gamma)
\ \Theta_\gamma'(t'),
\tag 8.10
$$
for all $\beta \in \N^{n-1}$,
after expanding with respect to $(w,z,\zeta,\xi)$
the power series appearing 
in $\underline{\Cal L}^\beta
(\bar g(\zeta, \xi)^\gamma)$, $\underline{\Cal L}^\beta
\bar f(\zeta, \xi)$ and after replacing $\xi$ by 
$z-i\Theta(\zeta,w,z)$, and similarly, we set\,:
$$
S_\beta(w,z,\zeta,t'):=
\underline{\Cal L}^\beta \bar f(\zeta, \xi) +i\sumg
{w'}^\gamma \ \underline{\Cal L}^\beta (\bar\Theta_\gamma'(\bar h(\zeta, \xi))),
\tag 8.11
$$
in coherence with the notation in eq.~\thetag{8.5} and finally also,
we set\,:
$$
T_\beta(w,z,\zeta,t'):=
-\underline{\omega}_\beta(w,\zeta,\xi)+\Theta_\beta'(t')+
\sumg {(\beta+\gamma)!\over \beta!\ \gamma!} \bar g(\zeta, \xi)^\gamma 
\ \Theta_{\beta+\gamma}'(t').
\tag 8.12
$$
We first remark that, by the very definition of 
$s_\beta$ and of $S_\beta$, we have
$$
{\partial s_\beta\over \partial t_j'} 
(w,\zeta, h(w,i\bar\Theta(w,\zeta,0)))\equiv
{\partial S_\beta\over \partial t_j'}
(w,i\bar\Theta(w,\zeta,0),\zeta, h(w,i\bar\Theta(w,\zeta,0))),
\tag 8.13
$$ 
as formal
power series, for all $\beta \in\N^{n-1}$, $1\leq j\leq n$.
Next, let us establish a useful correspondence between the vanishing
of the generic ranks of $(R_\beta)_{\v \beta\v \leq 2\kappa_0}$,
$(S_\beta)_{\v \beta\v \leq 2\kappa_0}$ and
$(T_\beta)_{\v \beta\v \leq 2\kappa_0}$.

\proclaim{Lemma~8.14}
The following properties are equivalent\,:
\roster
\item"{\bf (1)}" $\hbox{det}\left(
{\partial R_{\beta^i} \over \partial t_j'}
(w,z,\zeta,h(w,z))\right)_{1\leq i,j\leq n}
\equiv 0$, $\forall \ \beta^1,\ldots, \beta^n$,
$\vert \beta^1\vert,\ldots,\vert
\beta^n\vert \leq 2\kappa_0$.
\item"{\bf (2)}" $\hbox{det}\left(
{\partial S_{\beta^i} \over \partial t_j'}
(w,z,\zeta,h(w,z))\right)_{1\leq i,j\leq n}
\equiv 0$, $\forall \ \beta^1,\ldots, \beta^n$,
$\vert \beta^1\vert,\ldots,\vert
\beta^n\vert \leq 2\kappa_0$.
\item"{\bf (3)}" $\hbox{det}\left(
{\partial T_{\beta^i} \over \partial t_j'}
(w,z,\zeta,h(w,z))\right)_{1\leq i,j\leq n}
\equiv 0$, $\forall \ \beta^1,\ldots, \beta^n$,
$\vert \beta^1\vert,\ldots,\vert
\beta^n\vert \leq 2\kappa_0$.
\endroster
\endproclaim

\demo{End of proof of Lemma~8.8}
The proof of Lemma~8.14 will be given just below.  To finish the proof
of Lemma~8.8, we assume by contradiction that~\thetag{8.9} is untrue,
{\it i.e.} that {\bf (2)} of Lemma~8.14 holds with 
$z=i\bar\Theta(w,\zeta,0)$. 
According to {\bf (3)} of this lemma,
we also have that the generic rank of the $n\times
{(2\kappa_0+n-1)!\over (2\kappa_0)! (n-1)!}$ matrix
$$
\left\{
\aligned
\Cal N_{2\kappa_0}(w,\zeta):= 
&
\left(
{\partial \Theta_\beta' \over \partial t_j'}
(h(w,i\bar\Theta(w,\zeta,0)))+\right.\\ 
&
+\left. \sumg
{(\beta+\gamma)!\over \beta!\ \gamma!} \ \bar g(\zeta, 0)^\gamma \
{\partial \Theta_{\beta+\gamma}' \over \partial t_j'}
(h(w,i\bar\Theta(w,\zeta,0)))
\right)_{1\leq j\leq n}^{\vert \beta\vert \leq 2\kappa_0}
\endaligned\right.
\tag 8.15
$$
is strictly less than $n$. After making some obvious linear combinations
between the columns of $\Cal N_{2\kappa_0}$ with coefficients
being formal power series in $\zeta$ which are polynomial
with respect to the
$\bar g_j(\zeta,0)\in \frak{m}(\zeta)$, $1\leq j\leq n-1$, we can reduce 
$\Cal N_{2\kappa_0}$ to the matrix of same formal generic rank
$$
\left\{
\aligned
&
\Cal N_{2\kappa_0}^0(w,\zeta)
:=
\left(
{\partial \Theta_\beta' \over \partial t_j'}
(h(w,i\bar\Theta(w,\zeta,0)))+\right.\\ 
&
+\left. \sum_{\vert \gamma\vert \geq 2\kappa_0+1-\vert \beta\vert}
{(\beta+\gamma)!\over \beta!\ \gamma!} \ \bar g(\zeta, 0)^\gamma \
{\partial \Theta_{\beta+\gamma}' \over \partial t_j'}
(h(w,i\bar\Theta(w,\zeta,0)))
\right)_{1\leq j\leq n}^{\vert \beta\vert \leq 2\kappa_0}.
\endaligned\right.
\tag 8.16
$$
Now, taking the submatrix $\Cal N_{\kappa_0}^0$ of $\Cal
N_{2\kappa_0}^0$ for which $\v \beta \v \leq \kappa_0$, we see that we
have reduced $\Cal N_{\kappa_0}^0$ to the simpler form
$$
\Cal N_{\kappa_0}^0(w,\zeta)\equiv \Cal N^1_{\kappa_0}
(w,\zeta) \ \hbox{mod} \ (\frak{m} (\zeta)^{\kappa_0+1}
\Cal Mat_{n\times {(\kappa_0+n-1)!\over (\kappa_0)! (n-1)!}} 
(\C\dl w,\zeta\dr)),
\tag 8.17
$$
where
$$
\Cal N^1_{\kappa_0}(w,\zeta):=\left(
{\partial \Theta_\beta'\over \partial t_j'}
(h(w,i\bar\Theta(w,\zeta,0))\right)_{1\leq j\leq n}^{\vert \beta \vert
\leq \kappa_0}.
\tag 8.18
$$
But
$$ 
\hbox{det} \left( {\partial \Theta_{\underline{\beta}^i}'\over
\partial t_j'} (h(w,i\bar\Theta(w,\zeta,0)))
\right)_{1\leq i,j\leq n} \not\equiv 0 
\ \ \hbox{in} \ \C\dl w,\zeta\dr \ (\hbox{mod} 
\ \frak{m}(\zeta)^{\kappa_0+1}),
\tag 8.19
$$
by the choice of the 
$\underline{\beta}^i$'s and of $\kappa_0$,
which is the desired contradiction.
\qed
\enddemo
\enddemo

\demo{Proof of Lemma~8.14}
The equivalence {\bf (1)} $\Longleftrightarrow$ {\bf (3)} follows by
an inspection of the proof of Lemma~5.3\,: to pass from the system
$R_\beta=0$, $\vert
\beta \vert \leq 2\kappa_0$, to the system $T_\beta=0$, $\vert \beta
\vert \leq 2\kappa_0$, we have only use in the proof some
linear combinations with coefficients in $\C\dl \zeta,\xi\dr$.  The
equivalence {\bf (1)} $\Longleftrightarrow$ {\bf (2)} is related with
the substance of Lemma~5.10. Indeed, in the relation
$r'(t',\tau')\equiv
\alpha'(t',\tau') \, \bar r'(\tau',t')$, with $\alpha'(0,0)=-1$,
insert first $\tau':= \bar h(\tau)$ to get $r'(t',\bar h(\tau))\equiv
\alpha'(t',\bar h(\tau)) \, \bar r'(\bar h(\tau), t')$ and then
differentiate by the operator $\underline{\Cal L}^\beta$ to obtain
$$
\underline{\Cal L}^\beta r'(t',\bar h(\tau))\equiv 
\alpha'(t',\bar h(\tau)) \ \underline{\Cal L}^\beta 
\bar r'(\bar h(\tau), t')+
\sum_{\gamma\leq \beta, \gamma\neq\beta}
{\alpha'}_\gamma^\beta(t',t,\tau) \
\underline{\Cal L}^\gamma\bar r'(\bar h(\tau),t').
\tag 8.20
$$
In our notations, $R_\beta(w,z,\zeta,t')=\underline{\Cal L}^\beta
r'(t',\bar h(\tau))$ and
$S_\beta(w,z,\zeta,t')=\underline{\Cal L}^\beta \bar r'(\bar
h(\tau), t')$, after replacing $\xi$ by $z-i\Theta(\zeta,w,z)$.  We
deduce
$$
\left\{
\aligned
&
{\partial R_\beta\over \partial t_j'}(w,z,\zeta,h(w,z))
=\alpha'
(h(w,z),\bar h(\tau))
{\partial S_\beta\over \partial t_j'}(w,z,\zeta,h(w,z))+\\
&
+\sum_{\gamma\leq \beta, \gamma\neq\beta}{\alpha'}_\gamma^\beta(t,\tau,
h(w,z)) \ {\partial S_\gamma\over \partial t_j'}(w,z,\zeta,h(w,z)).
\endaligned\right.
\tag 8.21
$$
Equation~\thetag{8.21} shows that the terms ${\partial R_\beta\over
\partial t_j'}(w,z,\zeta,h(w,z))$ are trigonal linear combinations 
of the terms ${\partial S_\gamma\over \partial t_j'}(w,z,\zeta,h(w,z))$,
$\gamma\leq \beta$, with nonzero diagonal coefficients. This completes
the proofs of Lemmas~8.14 and~8.2 and completes finally our proof of
Theorem~1.2.
\qed
\enddemo
\enddemo

\remark{Remark}
Once we know that $h(w,z)\in \C\{w,z\}$, we deduce that
the reflection function associated with the formal equivalenve of
Theorem~1.2 is convergent, {\rm i.e.} that
$\Cal R_h'(w,z,\bar\lambda,\bar\mu)\in \C\{w,z,\bar\lambda,\bar\mu\}$.
\endremark

\enddocument
\end